\documentstyle[amsfonts,10pt]{article}
\input tcilatex
\begin{document}

\title{THE DOUBLE COVERS OF KLEIN SURFACES WITH NODES}
\author{IGNACIO C. GARIJO \\
Dpto. de Matem\'{a}ticas Fundamentales\\
Facultad de Ciencias, UNED\\
c/ Senda del Rey, 9\\
28040-Madrid. SPAIN}
\date{}
\maketitle

\begin{abstract}
{\small Given a Klein surface }$Y${\small , there is a unique symmetric
Riemann surface }$X${\small \ being the complex double of }$Y${\small . In
this article we shall show that the situation is not the same when we work
in the category of surfaces with nodes.}
\end{abstract}

Riemman surfaces with nodes appear as new objects of the Deligne-Mumford
compactification of the moduli spaces of smooth complex curves of given
genus \cite{D-M}. These compactifications are an useful tool in recent
important progress in mathematics. Concretely, the proof of Witten's
conjecture given by M. Kontsevich (see \cite{K}, \cite{L-Z} and \cite{Z})
uses, in an essential way, the cell descomposition of the spaces $\overline{%
{\cal M}}_{g,n}\times {\Bbb R}_{+}^n$, where $\overline{{\cal M}}_{g,n}$ is
the Deligne-Mumford compactification by Riemann surfaces with nodes of the
moduli space of smooth complex curves of genus $g$ with $n$ marked points.

Meanwhile, we are interested in the compactification of the moduli spaces of
smooth real algebraic curves of given genus $g$, given by \cite{Se} and \cite
{Si}. A smooth complex projective curve is simply a compact Riemann surface.
If the curve is defined by real polinomials, the corresponding Riemann
surface carries an antianalytic involution induced by the complex
conjugation and it is called a symmetric Riemann surface. Conversely, any
compact symmetric Riemann surface can be embedded in a complex projective
space in such a way that its image is a curve defined by real polinomials.
We conclude that a projective smooth real algebraic curve is simply a
compact symmetric Riemann surface.

Let be now $\left( X,\sigma \right) $ a symmetric Riemann surface where $X$
is a Riemann surface and $\sigma :X\rightarrow X$ is an antianalytic
involution. The quotient space $X/\left\langle \sigma \right\rangle $ has a
unique structure of Klein surface such that the projection $\pi
:X\rightarrow X/\left\langle \sigma \right\rangle $ is a morphism. On the
other hand, given a Klein surface $Y$, there exists a unique symmetric
Riemann surface $\left( X,\sigma \right) $ called the complex double of $Y$
satisfying that $X/\left\langle \sigma \right\rangle $ is isomorphic to $Y$.
By virtue of these observations, symmetric Riemann surfaces are the same
objects as Klein surfaces. The situation changes when we pass to the
category of surfaces with nodes. We shall prove that given a Klein surface
with nodes $Y$ there exist different symmetric Riemann surfaces with nodes $%
\left( X_i,\sigma _i\right) $ being the complex doubles of $Y$. Hence,
symmetric Riemann surfaces with nodes and Klein surfaces with nodes are
different objects.

The moduli spaces of smooth real algebraic curves of given genus $g$,
designed by ${\cal M}_{{\Bbb R}}^g$, are neither compact nor connected. In 
\cite{Se}, M. Sepp\"{a}l\"{a} constructs moduli spaces of stable symmetric
Riemann surfaces and proves that these spaces are compact and connected. In
other words, he constructs a compactification of ${\cal M}_{{\Bbb R}}^g$
using symmetric Riemann surfaces with nodes. If we construc the moduli space
of stable Klein surfaces then we can obtain another compactificaci\'{o}n of $%
{\cal M}_{{\Bbb R}}^g$.

We begin section 1 with the definition of Riemann surfaces with nodes. Each
point of this surfaces has a neigbourhood homeomorphic either to an open set
in the complex plane ${\Bbb C}$ or to 
\[
{\cal M}=\left\{ \left( z,w\right) \in {\Bbb C}^2\mid \quad z\cdot w=0,\quad
\left| z\right| <1,\quad \left| w\right| <1\right\} . 
\]
The points corresponding to the second case are the nodes of the surface.
The next step is constructing a normalization $\widehat{X}$ of the Riemann
surface with nodes $X$ ``unidentifying'' the nodes. The important property
of this normalization is that each continuous map $f:X_1\rightarrow X_2$
between Riemann surfaces with nodes induces a unique continuous map $%
\widehat{f}:\widehat{X_1}\rightarrow \widehat{X_2}$ between its
normalizations, called the lifting of $f$. Finally, a morphism $%
f:X_1\rightarrow X_2$ between Riemann surfaces with nodes is a continuous
map satisfying that

\begin{itemize}
\item  $f^{-1}\left( N\left( X_2\right) \right) \subset N\left( X_1\right) $
where $N\left( X_i\right) $ is the set of nodes of $X_i$,

\item  the lifting $\widehat{f}$ is a morphism,

\item  if $x\in N\left( X_1\right) $ and $f\left( x\right) \in N\left(
X_2\right) $, then the image by $f$ of a neigbourhood of $x$ is a
neigbourhood of $f\left( x\right) $.
\end{itemize}

In section 2 we define Klein surfaces with nodes and the objects associated
to them. Given an antianalytic involution $\sigma :X\rightarrow X$ on a
Riemann surface with nodes, the orbit space is a Klein surface with three
types of singularities (as quotients of the nodes): conic nodes (as in
Riemann surfaces with nodes), boundary nodes (half of a node in the
boundary) and inessential nodes (plane situations). The construction of
normalizations and liftings of maps and the definition of morphisms are
similar to those in section 1.

One of the main results appears in section 3. We begin with the definition
of a complex double of a Klein surface with nodes $Y$ as a triple $\left(
X,\pi ,\sigma \right) $ such that $\left( X,\sigma \right) $ is a symmetric
Riemann surface with nodes and $\pi :X\rightarrow Y$ is an unbranched double
cover satisfying that $\pi \circ \sigma =\pi $ and $\pi \left( N\left(
X\right) \right) =N\left( Y\right) $.

Now, if $z$ is a conic node of $Y$, then it generates two points, $z_1$ and $%
z_2$, in the normalization $\widehat{Y}$. Since $\widehat{Y}$ is a Klein
surface, then there exists its complex double $\left( \widehat{X},\widehat{%
\pi },\widehat{\sigma }\right) $. If $\widehat{\pi }^{-1}\left( z_i\right)
=\left\{ w_{i,1},w_{i,2}\right\} $, then we have two possible
identifications of each pair of points in order to construct a complex
double $\left( X,\pi ,\sigma \right) $. With these ideas in mind we can
construct $2^{\#N\left( Y,1\right) }$ non isomorphic complex doubles of $Y$
where $N\left( Y,1\right) $ are the set of conic nodes of $Y$, moreover,
these are the unique complex doubles. Then it does not exist unicity for the
complex doubles of a Klein surface with nodes $Y$.

We end this paper with three examples of Klein surfaces and its complex
doubles showing the variety of situations. In the first case the two complex
doubles are isomorphic as symmetric Riemann surfaces with nodes, in the
second example the two complex doubles are not homeomorphic as surfaces with
nodes with an orientation reversing involution and, in the last, the two
complex doubles are not homeomorphic as surfaces with nodes.

This article is part of my Ph. D. Thesis \cite{G} about Riemann and Klein
surfaces with nodes written under the supervision of Dr. Antonio F. Costa.

\section{Riemann surfaces with nodes.}

We begin this article defining the most important objects, notations and
constructions.

\begin{definition}
A surface with nodes, $\Sigma $, is a topological Haussdorff space such that
each point has a nieghbourhood homeomorphic either to an open set in the
complex plane ${\Bbb C}$ or to 
\[
{\cal M}=\left\{ \left( z,w\right) \in {\Bbb C}^2\mid \quad z\cdot w=0,\quad
\left| z\right| <1,\quad \left| w\right| <1\right\} . 
\]

We shall call a node of $\Sigma $ to a point such that each neighbourhood
contains an open set homeomorphic to ${\cal M}$. We denote the set of nodes
of $\Sigma $ by $N\left( \Sigma \right) $. Each connected component of $%
\Sigma \setminus N\left( \Sigma \right) $ is called a part of $\Sigma $.

A stable surface is a surface with nodes in which the Euler characteristic
of each part is negative.

Let $\Sigma $ be a compact, connected surface with $n$ nodes and Euler
characteristic $\chi \left( \Sigma \right) $. We define the genus of $\Sigma 
$ by: 
\[
g\left( \Sigma \right) =\left\{ 
\begin{array}{ll}
\frac 12\left( 2+n-\chi \left( \Sigma \right) \right) & \text{if }\Sigma 
\text{ is orientable,} \\ 
2+n-\chi \left( \Sigma \right) & \text{if }\Sigma \text{ is non-orientable.}
\end{array}
\right. 
\]
\end{definition}

With these definitions it is easy to show that $N\left( \Sigma \right) $ is
discrete and if $\Sigma $ is compact then both, the number of nodes and the
number of parts, is finite.

Let $\Sigma $ be a surface with nodes, then we can construct a unique
topological surface, $\Sigma _E$, substituting for each node a neighbourhood
homeomorphic to ${\cal M}$ by a cylinder. We can observe that $\chi \left(
\Sigma _E\right) =\chi \left( \Sigma \right) -n$ and hence $g\left( \Sigma
\right) =g\left( \Sigma _E\right) $.

The next step is to define an analytic structure on a surface with nodes:

\begin{definition}
Let $\Sigma $ be a surface with nodes; a chart is a pair $\left(
U_i,f_i\right) $ where $U_i$ is an open set of $\Sigma $ and $%
f_i:U_i\rightarrow V_i$ is an homeomorphism where $V_i$ is one of these sets:

\begin{description}
\item[I)]  An open set of ${\Bbb C}$.

\item[II)]  The disjoint union of two open sets of ${\Bbb C}$, $V_{i,1}$, $%
V_{i,2}$, with a point in each one that is identifyed, i.e. we take $z_i\in
V_{i,1}\cap V_{i,2}$ and we construct: 
\[
V_i=V_{i,1}\stackunder{z_i}{\sqcup }V_{i,2}=\left( V_{i,1}\times \left\{
1\right\} \cup V_{i,2}\times \left\{ 2\right\} \right) /\sim 
\]
where $\sim $ is the identification $\left( z_i,1\right) \sim \left(
z_i,2\right) $.
\end{description}

We say that two charts, $\left( U_i,f_i\right) $, $\left( U_j,f_j\right) $,
have analytical transition if $U_i\cap U_j=\emptyset $ or the transition map 
\[
f_{ij}=f_j\circ f_i^{-1}:f_i\left( U_i\cap U_j\right) \rightarrow f_j\left(
U_i\cap U_j\right) 
\]
is analytic in the image of the complementary of the nodes.

A Riemann surface with nodes is a pair $X=\left( \Sigma ,{\cal U}\right) $
wher $\Sigma $ is a surface with nodes and ${\cal U}$ is an analytical and
maximal atlas.

A stable Riemann surface is a Riemann surface with nodes such that $\Sigma $
is a stable surface.
\end{definition}

We can observe that a Riemann surface with nodes, $X$, is orientable and $%
X\setminus N\left( X\right) $ have structure of Riemann surface; indeed,
this definition extends the classical definition of Riemann surface.

\begin{remark}
Some authors call Riemann surface with nodes to a surface with nodes $\Sigma 
$ so that $\Sigma \setminus N\left( \Sigma \right) $ have an analytic
structure. This is a different definition. Let be $\Sigma ={\Bbb C}%
\stackunder{0}{\sqcup }{\Bbb C}$, $A={\Bbb C}\setminus \left\{ 0\right\} $
and $B=\left\{ z\in {\Bbb C}\mid \quad \left| z\right| >1\right\} $. Let be
now the following two charts: 
\[
\begin{array}{cccc}
f_1: & \Sigma \setminus N\left( \Sigma \right) & \rightarrow & A\times
\left\{ 1,2\right\} \\ 
& \left( z,i\right) & \rightarrow & \left( z,i\right) ,
\end{array}
\quad 
\begin{array}{cccc}
f_2: & \Sigma \setminus N\left( \Sigma \right) & \rightarrow & B\times
\left\{ 1,2\right\} \\ 
& \left( z,i\right) & \rightarrow & \left( z+\frac z{\left\| z\right\|
},i\right) .
\end{array}
\]
Then we have two non-equivalent analytic structures on $\Sigma \setminus
N\left( \Sigma \right) $ because the transition map 
\[
\begin{array}{cccc}
f_{12}=f_2\circ f_1^{-1}: & A\times \left\{ 1,2\right\} & \rightarrow & 
B\times \left\{ 1,2\right\} \\ 
& \left( z,i\right) & \rightarrow & \left( z+\frac z{\left\| z\right\|
},i\right)
\end{array}
\]
is not analytic. The first map can extend to an analytic structure on $%
\Sigma $ but it is impossible in the second case.

Studying Riemann surfaces with nodes we are interesting on the surfaces of
the compactification of the moduli space. Our definition describes the
surfaces appearing in the Deligne-Mumford compactification.
\end{remark}

Now we are going to construct a Riemann surface associated to each Riemann
surface with nodes that will play an important role in our work. Let $%
X=\left( \Sigma ,{\cal U}\right) $ be a Riemann surface with nodes and $%
N\left( \Sigma \right) =\left\{ z_i\right\} _{i\in I}$. Let $%
f_i:U_i\rightarrow V_i$ be the charts with $V_i=V_{i,1}\stackunder{w_i}{%
\sqcup }V_{i,2}$ and $U_i=U_{i,1}\stackunder{z_z}{\sqcup }U_{i,2}$ where $%
U_{i,k}=f_i^{-1}\left( V_{i,k}\right) $. We take $\Sigma \setminus N\left(
\Sigma \right) $ and construct $\widehat{X}=\left( \widehat{\Sigma },%
\widehat{{\cal U}}\right) $ identifying $U_{i,k}\backslash \left\{
z_i\right\} $ with $V_{i,k}$. Now we assign charts in the obvious way. Then
we have a Riemann surface $\widehat{X}$ such that, in general, is not
connected; in fact, there is a bijection between the connected components of 
$\widehat{X}$ and the parts of $X$. This new surface is called the
normalization of the Riemann surface with nodes $X$. The projection 
\[
\begin{array}{cccc}
p: & \widehat{X} & \rightarrow & X \\ 
& z & \rightarrow & z
\end{array}
\]
is a continuous, closed and onto map; in fact, it is an identification map.
Moreover $\#p^{-1}\left( z\right) =2$ if and only if $z\in N\left( \Sigma
\right) $. If $z\notin N\left( \Sigma \right) $ then $\#p^{-1}\left(
z\right) =1$, so $p$ has finite fibers and hence it is a perfect map. In
this situation we have that $X$ is compact if and only if $\widehat{X}$ is
compact.

It is easy to show that the structure defined in the previous remark do not
admit this normalization.

Using the normalization we can easily prove the following proposition:

\begin{proposition}
Each Riemann surface with nodes is obtained identifying pairs of points of a
discrete subset of a Riemann surface; moreover, each construction in this
way is a Riemann surface with nodes.
\end{proposition}

The following definition is for describing the maps between Riemann surfaces
with nodes.

\begin{definition}
Let $\Sigma _1,\Sigma _2$ be surfaces with nodes and $N_i=N\left( \Sigma
_i\right) $. A map between surfaces with nodes is a continuous map $f:\Sigma
_1\rightarrow \Sigma _2$ satisfying that $f^{-1}\left( N_2\right) \subset
N_1.$
\end{definition}

If we have a map $f:\Sigma _1\rightarrow \Sigma _2$ between surfaces with
nodes then we have the following diagram 
\[
\begin{array}{ccc}
\widehat{\Sigma _1} &  & \widehat{\Sigma _2} \\ 
\left\downarrow p_1\right. &  & \left\downarrow p_2\right. \\ 
\Sigma _1 & \stackrel{f}{\rightarrow } & \Sigma _2
\end{array}
\]
We shall define a map $\widehat{f}:\widehat{\Sigma _1}\rightarrow \widehat{%
\Sigma _2}$. Let $z$ be a point in $\widehat{\Sigma _1}$; if $f\circ
p_1\left( z\right) \notin N_2$, then $p_2^{-1}\circ f\circ p_1\left(
z\right) =\left\{ w\right\} $ and we define $\widehat{f}\left( z\right) =w$;
if $f\circ p_1\left( z\right) \in N_2$, then $p_1\left( z\right) \in N_1$, $%
p_1^{-1}\left( p_1\left( z\right) \right) =\left\{ z,z^{\prime }\right\} $
and $p_2^{-1}\circ f\circ p_1\left( z\right) =\left\{ w_1,w_2\right\} $ and
we define $\widehat{f}\left( z\right) =w_i$ in order to be continuous. We
call to this map the lifting of $f$ and we have that $f\circ p_1=p_2\circ 
\widehat{f}$. If there is other continuous map $\widehat{f^{\prime }}:%
\widehat{\Sigma _1}\rightarrow \widehat{\Sigma _2}$ such that $f\circ
p_1=p_2\circ \widehat{f^{\prime }}$ then $\widehat{f}\left( z\right) =%
\widehat{f^{\prime }}\left( z\right) \quad $for all $z\in \widehat{\Sigma _1}%
\setminus p_1^{-1}\left( N_1\right) $. Since $\widehat{\Sigma _1}\setminus
p_1^{-1}\left( N_1\right) $ is dense in $\widehat{\Sigma _1}$, we conclude
that $\widehat{f}=\widehat{f^{\prime }}$. We have the following result:

\begin{proposition}
Let $f:\Sigma _1\rightarrow \Sigma _2$ be a map between surfaces with nodes,
then there is a unique continuous map $\widehat{f}:\widehat{\Sigma _1}%
\rightarrow \widehat{\Sigma _2}$ such that $f\circ p_1=p_2\circ \widehat{f}$.
\end{proposition}

This lifting plays an important role in the theory of Riemann surfaces with
nodes because we can deduce properties of $\Sigma _i$ and $f$ looking at the
properties of $\widehat{\Sigma _i}$ and $\widehat{f}$. This lifting
satisfies that $\widehat{g\circ f}=\widehat{g}\circ \widehat{f}$, $\widehat{%
Id_\Sigma }=Id_{\widehat{\Sigma }}$, $f$ is bijective if and only $\widehat{f%
}$ is bijective and $\widehat{f^{-1}}=\widehat{f}^{-1}$.

Let now $z$ be in $\Sigma _1$. Then $\widehat{f}\left( p_1^{-1}\left(
z\right) \right) \subset p_2^{-1}\left( f\left( z\right) \right) $ and $%
1\leq \#\widehat{f}\left( p_1^{-1}\left( z\right) \right) \leq
\#p_2^{-1}\left( f\left( z\right) \right) \leq 2$. If $f\left( z\right)
\notin N_2$ then $\widehat{f}\left( p_1^{-1}\left( z\right) \right)
=p_2^{-1}\left( f\left( z\right) \right) $, but if $f\left( z\right) \in N_2$
then $\#p_2^{-1}\left( f\left( z\right) \right) =2$ and $z\in N_1$ and it
may occur that $\#\widehat{f}\left( p_1^{-1}\left( z\right) \right) =1$. In
this case $\widehat{f}$ does not ``complete'' $p_2^{-1}\left( N_2\right) $.
This situation motivates the following definition:

\begin{definition}
We say that a map $f:\Sigma _1\rightarrow \Sigma _2$ between surfaces with
nodes is complete if $\widehat{f}\left( p_1^{-1}\left( z\right) \right)
=p_2^{-1}\left( f\left( z\right) \right) .$
\end{definition}

Later we shall see that if $f$ is not complete then can happen strange
situations. Now we define the morphisms between Riemann surfaces with nodes:

\begin{definition}
We say that a map $f:X_1\rightarrow X_2$ between Riemann surfaces with nodes
is analytic (resp. antianalytic) if $f:X_1\setminus N_1\rightarrow
X_2\setminus N_2$ is an analytic (resp. antianalytic) morphism between
Riemann surfaces.

A morphism (resp. antianalytic morphism) between Riemann surfaces with nodes
is an analytic (resp. antianalytic) and complete map $f:X_1\rightarrow X_2$
between Riemann surfaces with nodes.
\end{definition}

We define isomorphisms and automorphisms in the natural way. With these
definitions we have interesting properties: $f$ is analytic (resp. morphism,
antianalytic, isomorphism,...) if and only if $\widehat{f}$ is analytic
(resp. morphism, antianalytic, isomorphism,...). Moreover, as in the theory
of Riemann surfaces, a morphism $f:X_1\rightarrow X_2$ that it is not
constant in each part of $X_1$ is an open map.

Now we show with an example that the property ``complete'' is essential. Let
be $X=\widehat{{\Bbb C}}\stackunder{0}{\sqcup }\widehat{{\Bbb C}}$ where $%
\widehat{{\Bbb C}}={\Bbb C}\cup \left\{ \infty \right\} $ is the Riemann
sphere and consider the map 
\[
\begin{array}{cccc}
f: & X & \rightarrow & X \\ 
& \left( z,i\right) & \rightarrow & \left( z,1\right) .
\end{array}
\]
This is an analytic map between Riemann surfaces with nodes but $f$ is not
complete. In this case $X$ is an open set but $f\left( X\right) =\widehat{%
{\Bbb C}}\times \left\{ 1\right\} $ is not, hence $f$ is not an open map.

\section{Klein surfaces with nodes.}

Klein surfaces appear as quotients of Riemann surfaces by antianalytic
involutions. Then we shall study the quotient of a Riemann surface with
nodes by an antianalytic involution in the neigbourhood of a node.

Let $X$ be a Riemann surface with nodes, $\sigma :X\rightarrow X$ an
antianalytic involution, $z\in N\left( X\right) $ and $V_1,V_2\subset X$ two
disks such that $U=V_1\stackunder{z}{\sqcup }V_2$ is a neigbourhood of $z$.
We are interested in the quotient $U/\left\langle \sigma \right\rangle $. It
can happen three situations:

\begin{itemize}
\item  $\sigma \left( z\right) \neq z$. Then $\sigma \left( z\right) \in
N\left( X\right) \setminus \left\{ z\right\} $ and the quotient $%
U/\left\langle \sigma \right\rangle \simeq U$, so $z$ is again a node.

\item  $\sigma \left( V_1\right) =V_2$. In this case $U/\left\langle \sigma
\right\rangle \simeq V_1$ and the quotient is plane.

\item  $\sigma \left( V_i\right) =V_i$. Then $\sigma \mid _{V_i}$ is a
disk's symmetry along a diameter containing $z$. Hence $W_i=V_i/\left\langle
\sigma \right\rangle $ is a half-disk, so $U/\left\langle \sigma
\right\rangle \simeq W_1\stackunder{z}{\sqcup }W_2$ and it is half of a node.
\end{itemize}

Using these ideas we are going to extend the concept of surface with nodes:

\begin{definition}
A surface with nodes is a pair $S=\left( \Sigma ,{\cal D}\right) $ where $%
\Sigma $ is a topological Haussdorff space and ${\cal D}\subset \Sigma $ is
a discrete set of distinguished points of $\Sigma $ satisfying that each
point of $\Sigma $ have a neigbourhood $U$ homeomorphic to one of the
following sets:

\begin{description}
\item[1)]  An open set of ${\Bbb C}$,

\item[2)]  An open set of ${\Bbb C}^{+}=\left\{ z\in {\Bbb C}\mid \quad 
\func{Im}\left( z\right) \geq 0\right\} $,

\item[3)]  ${\cal M}=\left\{ \left( z,w\right) \in {\Bbb C}^2\mid \quad
z\cdot w=0,\quad \left| z\right| <1,\quad \left| w\right| <1\right\} $,

\item[4)]  ${\cal M}^{+}=\left\{ \left( z,w\right) \in \left( {\Bbb C}%
^{+}\right) ^2\mid \quad z\cdot w=0,\quad \left| z\right| <1,\quad \left|
w\right| <1\right\} $,
\end{description}

and if $z\in {\cal D}$ then $U$ is homeomorphic to an open set of ${\Bbb C} $%
.

If $\left\{ \left( U_i,\varphi _i\right) \right\} _{i\in I}$ is an atlas of $%
\Sigma $ then we define the boundary of $\Sigma $ as 
\[
\partial \Sigma =\left\{ z\in \Sigma \mid \quad \text{there is }i\in I\text{
with }z\in U_i\text{, }\varphi _i\left( U_i\right) \subset {\Bbb C}^{+}\text{
and }\varphi _i\left( z\right) \in {\Bbb R}=\partial {\Bbb C}^{+}\right\}
\cup 
\]
\[
\cup \left\{ z\in \Sigma \mid \quad \text{there is }i\in I\text{ with }z\in
U_i\text{, }\varphi _i\left( U_i\right) \subset \left( {\Bbb C}^{+}\right) ^2%
\text{ and }\varphi _i\left( z\right) \in \partial \left( \left( {\Bbb C}%
^{+}\right) ^2\right) \right\} . 
\]

If $z\in \Sigma $ and there is $i\in I$ with $z\in U_i$, $\varphi _i\left(
U_i\right) ={\cal M}$ and $\varphi _i\left( z\right) =\left( 0,0\right) $
then we say that $z$ is a conic node or an 1-node. We denote the set of
conic nodes by $N\left( \Sigma ,1\right) $.

If $z\in {\cal D}$ then we say that $z$ is an inessential node or a 2-node.
We denote the set of inessential nodes by $N\left( \Sigma ,2\right) $.

If $z\in \Sigma $ and there is $i\in I$ with $z\in U_i$, $\varphi _i\left(
U_i\right) ={\cal M}^{+}$ and $\varphi _i\left( z\right) =\left( 0,0\right) $
then we say that $z$ is a boundary node or a 3-node. We denote the set of
boundary nodes by $N\left( \Sigma ,3\right) $.

If $z$ belongs to $N\left( \Sigma \right) =N\left( \Sigma ,1\right) \cup
N\left( \Sigma ,2\right) \cup N\left( \Sigma ,3\right) $, then we say that $%
z $ is a node.

Finally, we shall call a part of $\Sigma $ to each connected component of $%
\Sigma \setminus N\left( \Sigma \right) $.
\end{definition}

As in the original definition $N\left( \Sigma \right) $ is discrete, ${\cal D%
}\cap \partial \Sigma =\emptyset $ and if $\Sigma $ is compact then both,
the number of nodes and the number of parts, is finite.

Now we are going to define a dianalytic structure on a surface with nodes:

\begin{definition}
Let $\Sigma $ be a surface with nodes; a chart is a pair $\left(
U_i,f_i\right) $ where $U_i$ is an open set of $\Sigma $ and $%
f_i:U_i\rightarrow V_i$ is an homeomorphism wher $V_i$ is one of the
following sets:

\begin{description}
\item[I)]  An open set of ${\Bbb C}^{+}$.

\item[II)]  The disjoint union of two open sets of ${\Bbb C}^{+}$, $V_{i,1}$%
, $V_{i,2}$, identified across $z_i\in V_{i,1}\cap V_{i,2}$ and denoted by: 
\[
V_i=V_{i,1}\stackunder{z_i}{\sqcup }V_{i,2}=\left( V_{i,1}\times \left\{
1\right\} \cup V_{i,2}\times \left\{ 2\right\} \right) /\sim 
\]
where $\sim $ is the identification $\left( z_i,1\right) \sim \left(
z_i,2\right) $.
\end{description}

We say that two charts, $\left( U_i,f_i\right) $, $\left( U_j,f_j\right) $,
have dianalytical transition if $U_i\cap U_j=\emptyset $ or the transition
map 
\[
f_{ij}=f_j\circ f_i^{-1}:f_i\left( U_i\cap U_j\right) \rightarrow f_j\left(
U_i\cap U_j\right) 
\]
is dianalytic in the image of the complementary of the nodes.

A Klein surface with nodes is a triple $X=\left( \Sigma ,{\cal D},{\cal U}%
\right) $ wher $\left( \Sigma ,{\cal D}\right) $ is a surface with nodes and 
${\cal U}$ is a dianalytical and maximal atlas.
\end{definition}

This definition of Klein surface with nodes extends the definition of Klein
surface. Moreover, each Riemann surface with nodes define a unique structure
of Klein surface with nodes, and each orientable Klein surfaces with nodes $%
\left( \Sigma ,\emptyset ,{\cal U}\right) $ without boundary and without
inessential nodes defines two structures of Riemann surface with nodes $%
\left( \Sigma ,{\cal U}_1\right) $ and $\left( \Sigma ,\overline{{\cal U}_1}%
\right) $. In this sense, we say that the definition of Klein surface with
nodes extends the definition of Riemann surface with nodes.

As in the case of Riemann surfaces with nodes, if we have a Klein surface
with nodes $X=\left( \Sigma ,{\cal D},{\cal U}\right) $ then we can
construct its normalization $\widehat{X}=\left( \widehat{\Sigma },\widehat{U}%
\right) $ using a similar way as those described in section 1. In this case
the projection 
\[
\begin{array}{cccc}
p: & \widehat{X} & \rightarrow & X \\ 
& z & \rightarrow & z
\end{array}
\]
have the same properties that in the case of Riemann surfaces with nodes. In
particular, $\#p^{-1}\left( z\right) =2$ if and only if $z\in N\left( \Sigma
,1\right) \cup N\left( \Sigma ,3\right) $, in other case $\#p^{-1}\left(
z\right) =1$.

We have the following result:

\begin{proposition}
Each Klein surface with nodes is obtained by a process of identification of
pairs of points and distinction in a discrete subset of points of a Klein
surface following the next rules:

\begin{description}
\item[1.]  The distinguished points can not belong to the boundary of the
Klein surface,

\item[2.]  We can not identify a point belonging to the boundary with a
point not belonging to the boundary.
\end{description}

Moreover, each topological space constructed in this way is a Klein surface
with nodes.
\end{proposition}

Now we shall deal with maps.

\begin{definition}
Let $\Sigma _1,\Sigma _2$ be surfaces with nodes. A map between surfaces
with nodes is a continuous map $f:\Sigma _1\rightarrow \Sigma _2$ such that $%
f^{-1}\left( N\left( \Sigma _2,i\right) \right) \subset N\left( \Sigma
_1,i\right) \cup N\left( \Sigma _1,1\right) $ and $f\left( \partial \Sigma
_1\right) \subset \partial \Sigma _1$.
\end{definition}

As in the previous section, we have the following proposition:

\begin{proposition}
Let $f:\Sigma _1\rightarrow \Sigma _2$ be a map between surfaces with nodes,
then there is a unique continuous map $\widehat{f}:\widehat{\Sigma _1}%
\rightarrow \widehat{\Sigma _2}$ such that $f\circ p_1=p_2\circ \widehat{f}$.
\end{proposition}

As in section 1, if $z\in \Sigma _1$, then $\widehat{f}\left( p_1^{-1}\left(
z\right) \right) \subset p_2^{-1}\left( f\left( z\right) \right) $ and $%
1\leq \#\widehat{f}\left( p_1^{-1}\left( z\right) \right) \leq
\#p_2^{-1}\left( f\left( z\right) \right) \leq 2$. If $f\left( z\right)
\notin N\left( \Sigma _2,1\right) \cup N\left( \Sigma _2,3\right) $ then $%
\widehat{f}\left( p_1^{-1}\left( z\right) \right) =p_2^{-1}\left( f\left(
z\right) \right) $, but if $f\left( z\right) \in N\left( \Sigma _2,1\right)
\cup N\left( \Sigma _2,3\right) $ then $\#p_2^{-1}\left( f\left( z\right)
\right) =2$. In this case $z\in N\left( \Sigma _1,1\right) \cup N\left(
\Sigma _1,3\right) $ and may occur that $\#\widehat{f}\left( p_1^{-1}\left(
z\right) \right) =1$. Hence we have the definitions:

\begin{definition}
We say that a map $f:\Sigma _1\rightarrow \Sigma _1$ between surfaces with
nodes is complete if $\widehat{f}\left( p_1^{-1}\left( z\right) \right)
=p_2^{-1}\left( f\left( z\right) \right) $.

We say that a map $f:X_1\rightarrow X_2$ between Klein surfaces with nodes
is dianalytic if $f:X_1\setminus N_1\rightarrow X_2\setminus N_2$ is a
dianalytic morphism between Klein surfaces.

A morphism between Klein surfaces with nodes is a dianalytic and complete
map $f:X_1\rightarrow X_2$ between Klein surfaces with nodes.
\end{definition}

We define isomorphisms and automorphisms in the natural way. With these
definitions we have that the maps $f$ and $\widehat{f}$ have the same
properties than in the case of Riemann surfaces with nodes.

\section{Symmetric Riemann surfaces with nodes and double covers.}

We begin this section with some definitions:

\begin{definition}
A symmetric Riemann surface with nodes is a pair $\left( X,\sigma \right) $
with $X$ a Riemann surface with nodes and $\sigma :X\rightarrow X$ an
antianalytic involution.

A map (resp. homeomorphism, morphism, antianalytic morphism,...) $f:\left(
X_1,\sigma _1\right) \rightarrow \left( X_2,\sigma _2\right) $ between
symmetric Riemann surfaces with nodes is a map (resp. homeomorphism,
morphism, antianalytic morphism,...) $f:X_1\rightarrow X_2$ such that $%
f\circ \sigma _1=\sigma _2\circ f$.
\end{definition}

If we have a symmetric Riemann surface with nodes $\left( X,\sigma \right) $
then $\left( \widehat{X},\widehat{\sigma }\right) $ is a symmetric Riemann
surface where $\widehat{X}$ is the normalization of $X$ and $\widehat{\sigma 
}$ is the lifting of $\sigma $. And, if $f:\left( X_1,\sigma _1\right)
\rightarrow \left( X_2,\sigma _2\right) $ is a map (resp. morphism,...)
between symmetric Riemann surfaces with nodes then, its lifting $\widehat{f}$
is a map (resp. morphism,...) between symmetric Riemann surfaces.

Let us remember a known theorem about complex doubles of Klein surfaces:

\begin{theorem}
\cite{A-G}. Let $X$ be a Klein surface, then there exists a triple $\left(
X_c,\pi _c,\sigma _c\right) $ such that $\left( X_c,\sigma _c\right) $ is a
symmetric Riemann surface and $\pi _c:X_c\rightarrow X$ is an unbranched
double cover satisfying that $\pi _c\circ \sigma _c=\pi _c$. If $\left(
X_c^{\prime },\pi _c^{\prime },\sigma _c^{\prime }\right) $ is other triple
with the same property, then there exists a unique analytic isomorphism $%
f:\left( X_c^{\prime },\sigma _c^{\prime }\right) \rightarrow \left(
X_c,\sigma _c\right) $ between symmetric Riemann surfaces such that $\pi
_c^{\prime }=\pi _c\circ f$. This triple, unique up isomorphism, is called
the complex double of $X$.
\end{theorem}

We are going two show that it does not happen the same when we consider
Klein surfaces with nodes.

\begin{definition}
We say that a morphism $f:X_1\rightarrow X_2$ is a double cover of Klein
surfaces with nodes if its lifting $\widehat{f}:\widehat{X_1}\rightarrow 
\widehat{X_2}$ is a double cover of Klein surfaces. The double cover $f$ is
branched or not depending on $\widehat{f}$.

Let $Y$ be a Klein surface with nodes. We say that the triple $\left( X,\pi
,\sigma \right) $ is a complex double of $Y$ if $\left( X,\sigma \right) $
is a symmetric Riemann surface with nodes and $\pi :X\rightarrow Y$ is an
unbranched double cover satisfying that $\pi \circ \sigma =\pi $ and $\pi
\left( N\left( X\right) \right) =N\left( Y\right) $.

We say that two complex doubles, $\left( X_1,\pi _1,\sigma _1\right) $, $%
\left( X_2,\pi _2,\sigma _2\right) $, of $Y$ are isomorphic if there is an
isomorphism $f:\left( X_1,\sigma _1\right) \rightarrow \left( X_2,\sigma
_2\right) $ between symmetric Riemann surfaces with nodes such that $\pi
_2\circ f=\pi _1$.
\end{definition}

We shall show with an example the importance of the condition $\pi \left(
N\left( X\right) \right) =N\left( Y\right) $. Let $\widehat{{\Bbb C}^{+}}=%
{\Bbb C}^{+}\cup \left\{ \infty \right\} $ be the disk, $\widehat{X}=%
\widehat{{\Bbb C}^{+}}\times \left\{ 1\right\} \cup \widehat{{\Bbb C}^{+}}%
\times \left\{ 2\right\} $ and $X=\widehat{X}/\sim $ where $\sim $ is the
identification $\left( 0,1\right) \sim \left( 0,2\right) $. Then $X$ is a
Klein surface with nodes and $\widehat{X}$ is its normalization. Let be now $%
\widehat{X_c}=\widehat{{\Bbb C}}\times \left\{ 1\right\} \cup \widehat{{\Bbb %
C}}\times \left\{ 2\right\} $ and we define: 
\[
\begin{array}{cccc}
\widehat{\sigma _c}: & \widehat{X_c} & \rightarrow & \widehat{X_c} \\ 
& \left( z,i\right) & \rightarrow & \left( \overline{z},i\right) ,
\end{array}
\qquad 
\begin{array}{cccc}
\widehat{\pi _c}: & \widehat{X_c} & \rightarrow & \widehat{X} \\ 
& \left( z,i\right) & \rightarrow & \left( \phi \left( z\right) ,i\right)
\end{array}
\]
where $\phi \left( z\right) =\func{Re}z+\left| \func{Im}z\right| \sqrt{-1}$.
Then $\left( \widehat{X_c},\widehat{\pi _c},\widehat{\sigma _c}\right) $ is
the complex double of $\widehat{X}$.

Let be now $X_1=\widehat{X_c}/\sim _1$ where $\sim _1$ is the identification 
$\left( 0,1\right) \sim _1\left( 0,2\right) $ and $X_2=X_1/\sim _2$ where $%
\sim _2$ is the identification $\left( i,1\right) \sim _2\left( -i,1\right) $%
. Now we define $\pi _i$ and $\sigma _i$ in the obvious way. Hence $\left(
X_i,\sigma _i\right) $, for $i=1,2$, are two symmetric Riemann surfaces with
nodes such that $\pi _i$ are unbranched double covers such that $\pi _i\circ
\sigma _i=\pi _i$, but $\pi _1\left( N\left( X_1\right) \right) =N\left(
X\right) $ and $\pi _2\left( N\left( X_2\right) \right) \neq N\left(
X\right) $. Its obvious that we can continue creating nodes on $X_1$
identifying $z$ with $\sigma _1\left( z\right) $.

Now we are going to see that if $Y$ is a Klein surface with nodes then, in
general, there is not a unique complex double $\left( X,\pi ,\sigma \right) $%
. Previously, let us see in what conditions a map between normalizations is
a lifting.

\begin{proposition}
Let $X_1$, $X_2$ be Klein surfaces with nodes, $\widehat{N}\left(
X_i,j\right) =p_i^{-1}\left( N\left( X_i,j\right) \right) $ and $g:\widehat{%
X_1}\rightarrow \widehat{X_2}$ be a continuous map such that $g^{-1}\left( 
\widehat{N}\left( X_2,j\right) \right) \subset \widehat{N}\left(
X_1,j\right) \cup \widehat{N}\left( X_1,1\right) $, $g\left( \partial 
\widehat{X_1}\right) \subset \partial \widehat{X_2}$ and $\#p_2\circ g\circ
p_1^{-1}\left( z\right) =1\quad $for all $z\in X_1$. With this conditions
there exists a unique $f:X_1\rightarrow X_2$ map between Klein surfaces with
nodes such that $g=\widehat{f}$.
\end{proposition}

\TeXButton{Proof}{\proof  }Since $\#p_2\circ g\circ p_1^{-1}\left( z\right)
=1$, then $g$ is compatible with the identification maps $p_1$ and $p_2$ and
exists a continuous map $f:X_1\rightarrow X_2$ such that $p_2\circ g=f\circ
p_1$. Now it is easy to see that $f\left( \partial X_1\right) \subset
\partial X_2$ and $f^{-1}\left( N\left( X_2,j\right) \right) \subset N\left(
X_1,j\right) \cup N\left( X_1,1\right) $. Hence $f$ is a map between Klein
surfaces with nodes and, for the unicity of the liftings, we have that $%
\widehat{f}=g$.

If there are $f_{1},f_{2}:X_{1}\rightarrow X_{2}$ maps between Klein
surfaces with nodes such that $\widehat{f_{1}}=\widehat{f_{2}}=g$, then $%
f_{1}\left( z\right) =f_{2}\left( z\right) \quad $for all $z\in
X_{1}\setminus N\left( X_{1}\right) $. But $X_{1}\setminus N\left(
X_{1}\right) $ is dense in $X_{1}$, hence $f_{1}=f_{2}$.\TeXButton{End Proof}
{\endproof  }

This proposition is also valid if we replace Klein surfaces with nodes by
Riemann surfaces with nodes.

Let now $X$ be a Klein surface with nodes, $\widehat{X}$ be its
normalization and $\left( \widehat{X_c},\widehat{\pi _c},\widehat{\sigma _c}%
\right) $ be its complex double. We define $X_c=\widehat{X_c}/\sim $ where $%
\sim $ are the identifications:

\begin{itemize}
\item  If $z\in N\left( X,3\right) $, then 
\[
p^{-1}\left( z\right) =\left\{ z_1,z_2\right\} \subset \partial \widehat{X}%
\text{ and }\widehat{\pi _c}^{-1}\left( z_i\right) =\left\{ w_i\right\} . 
\]
In this case $w_1\sim w_2$.

\item  If $z\in N\left( X,2\right) $, then 
\[
p^{-1}\left( z\right) =\left\{ z_1\right\} \notin \partial \widehat{X}\text{
and }\widehat{\pi _c}^{-1}\left( z_1\right) =\left\{ w_1,w_2\right\} \text{.}
\]
In this case $w_1\sim w_2$.

\item  If $z\in N\left( X,1\right) $, then 
\[
p^{-1}\left( z\right) =\left\{ z_1,z_2\right\} \notin \partial \widehat{X}%
\text{ and }\widehat{\pi _c}^{-1}\left( z_i\right) =\left\{
w_{i,1},w_{i,2}\right\} \text{.} 
\]
In this case there are two possible indentifications: $w_{1,1}\sim w_{2,1}$
and $w_{1,2}\sim w_{2,2}$ or $w_{1,1}\sim w_{2,2}$ and $w_{1,2}\sim w_{2,1}$.
\end{itemize}

We have constructed $2^{\#N\left( X,1\right) }$ Riemann surfaces with nodes
and the normalization of each of them is $\widehat{X_c}$.

Now we have the followings diagrams: 
\[
\begin{array}{ccc}
\widehat{X_c} & \stackrel{\widehat{\pi _c}}{\rightarrow } & \widehat{X} \\ 
\left\downarrow p_c\right. &  & \left\downarrow p\right. \\ 
X_c &  & X
\end{array}
\qquad \qquad 
\begin{array}{ccc}
\widehat{X_c} & \stackrel{\widehat{\sigma _c}}{\rightarrow } & \widehat{X}
\\ 
\left\downarrow p_c\right. &  & \left\downarrow p_c\right. \\ 
X_c &  & X_c
\end{array}
\]
Using the previous proposition and looking at the construction, it is easy
to show that there exist exactly two maps between Klein surfaces with nodes $%
\pi _c:X_c\rightarrow X$ and $\sigma _c:X_c\rightarrow X_c$ such that $%
\widehat{\pi _c}$ and $\widehat{\sigma _c}$ are their liftings. As $\widehat{%
\pi _c}$ is an unbranched double cover morphism and $\widehat{\sigma _c}$ is
an antianalytic morphism, so are $\pi _c$ and $\sigma _c$.

As $\widehat{\pi _c}\circ \widehat{\sigma _c}=\widehat{\pi _c}$ is the
lifting of $\pi _c\circ \sigma _c$ and $\pi _c$ then $\pi _c\circ \sigma
_c=\pi _c$; as $\widehat{\sigma _c}^2=Id_{\widehat{X_c}}$ is the lifting of $%
\sigma _c^2$ and $Id_{X_c}$ then $\sigma _c^2=Id_{X_c}$. Finally it is easy
to show that $\pi _c\left( N\left( X_c\right) \right) =N\left( X\right) $.
Then we have constructed $2^{\#N\left( X,1\right) }$ triples, $\left(
X_c,\pi _c,\sigma _c\right) $, that are complex doubles of $X$.

\begin{theorem}
Let $X$ be a Klein surface with nodes and $\left( Y,\pi ,\sigma \right) $ be
a complex double of $X$. Then there exists a unique $\left( X_c,\pi
_c,\sigma _c\right) $ complex double of $X$ constructed previously and there
exists a unique $f:\left( Y,\pi ,\sigma \right) \rightarrow \left( X_c,\pi
_c,\sigma _c\right) $ isomorphism between double covers.
\end{theorem}

\TeXButton{Proof}{\proof  }With these conditions there exists $\widehat{f}%
:\left( \widehat{Y},\widehat{\sigma }\right) \rightarrow \left( \widehat{X_c}%
,\widehat{\sigma _c}\right) $ isomorphism between symmetric Riemann surfaces
such that $\widehat{\pi _c}\circ \widehat{f}=\widehat{\pi }$. Then we have
the following commutative diagram: 
\[
\begin{array}{ccccc}
\widehat{Y} & \stackrel{\widehat{f}}{\rightarrow } & \widehat{X_c} & 
\stackrel{\widehat{\pi _c}}{\rightarrow } & \widehat{X} \\ 
p_y &  &  &  & \left\downarrow p_x\right. \\ 
Y &  & \stackrel{\pi }{\rightarrow } &  & X
\end{array}
\]
We construct $X_c=\widehat{X_c}/\sim _c$ where $\sim _c$ are the
identifications: 
\[
a\sim _cb\text{ if and only if }p_y\left( \widehat{f}^{-1}\left( a\right)
\right) =p_y\left( \widehat{f}^{-1}\left( b\right) \right) 
\]
and we call $p_c:\widehat{X_c}\rightarrow X_c$ to the projection. We are
going to see that $X_c$ is one of the $2^{\#N\left( X,1\right) }$ Riemann
surfaces with nodes constructed previously.

If $a_1\sim _ca_2$ then $p_x\circ \widehat{\pi _c}\left( a_i\right) =\pi
\circ p_y\circ \widehat{f}^{-1}\left( a_i\right) \in \pi \left( N\left(
Y\right) \right) =N\left( X\right) $. Now, if $a\in \widehat{\pi _c}%
^{-1}\circ p_x^{-1}\left( N\left( X\right) \right) $ then $\pi \circ
p_y\circ \widehat{f}^{-1}\left( a\right) \in N\left( X\right) $ and $%
p_y\circ \widehat{f}^{-1}\left( a\right) \in N\left( Y\right) $; hence there
is a unique $b\in \widehat{X_c}\backslash \left\{ a\right\} $ such that $%
p_y\circ \widehat{f}^{-1}\left( a\right) =p_y\circ \widehat{f}^{-1}\left(
b\right) $ and then 
\[
p_x\circ \widehat{\pi _c}\left( b\right) =\pi \circ p_y\circ \widehat{f}%
^{-1}\left( b\right) =\pi \circ p_y\circ \widehat{f}^{-1}\left( a\right)
=p_x\circ \widehat{\pi _c}\left( a\right) \in N\left( X\right) \text{.} 
\]

Hence we are going to study only $\widehat{\pi _c}^{-1}\circ p_x^{-1}\left(
N\left( X\right) \right) $.

\begin{itemize}
\item  If $z\in N\left( X,3\right) $ then 
\[
p_x^{-1}\left( z\right) =\left\{ z_1,z_2\right\} \subset \partial \widehat{X}%
\text{ and }\widehat{\pi _c}^{-1}\left( z_i\right) =\left\{ w_i\right\} 
\text{.} 
\]
Thus $\left( \widehat{\pi _c}^{-1}\circ p_x^{-1}\right) \left( z\right)
=\left\{ w_1,w_2\right\} $ and $w_1\sim _cw_2$.

\item  If $z\in N\left( X,2\right) $ then 
\[
p_x^{-1}\left( z\right) =\left\{ z_1\right\} \notin \partial \widehat{X}%
\text{ and }\widehat{\pi _c}^{-1}\left( z_1\right) =\left\{ w_1,w_2\right\}
. 
\]
Thus $\left( \widehat{\pi _c}^{-1}\circ p_x^{-1}\right) \left( z\right)
=\left\{ w_1,w_2\right\} $ and $w_1\sim _cw_2$.

\item  If $z\in N\left( X,1\right) $ then 
\[
p_x^{-1}\left( z\right) =\left\{ z_1,z_2\right\} \notin \partial \widehat{X}%
\text{ and }\widehat{\pi _c}^{-1}\left( z_i\right) =\left\{
w_{i,1},w_{i,2}\right\} . 
\]
Thus $\left( \widehat{\pi _c}^{-1}\circ p_x^{-1}\right) \left( z\right)
=\left\{ w_{i,j}\right\} _{i,j=1}^2$. If $p_y\left( \widehat{f}^{-1}\left(
w_{i,1}\right) \right) =p_y\left( \widehat{f}^{-1}\left( w_{i,2}\right)
\right) $ then $\pi $ is not complete because $\widehat{\pi }\left( \widehat{%
f}^{-1}\left( w_{i,1}\right) \right) =\widehat{\pi }\left( \widehat{f}%
^{-1}\left( w_{i,2}\right) \right) =z_i$. Hence there are two possible
identifications as in the previous construction.
\end{itemize}

Then $X_c$ is one of the $2^{\#N\left( X,1\right) }$ Riemann surfaces with
nodes constructed previously. As $\#p_c\circ f\circ p_y^{-1}\left( x\right)
=1\quad $for all $x\in Y$ then there is a unique $f:Y\rightarrow X_c$ map
between Riemann surfaces with nodes such that $\widehat{f}$ is its lifting.
Since $\widehat{f}$ is an isomorphism then $f$ is also an isomorphism.

We define $\sigma _{c}=f\circ \sigma \circ f^{-1}:X_{c}\rightarrow X_{c}$
and $\pi _{c}=\pi \circ f^{-1}:X_{c}\rightarrow X$. Then $\left( X_{c},\pi
_{c},\sigma _{c}\right) $ is one of the previous complex doubles and $%
f:\left( Y,\pi ,\sigma \right) \rightarrow \left( X_{c},\pi _{c},\sigma
_{c}\right) $ is an isomorphism between double covers.

Let $\left( X_1,\pi _1,\sigma _1\right) $, $\left( X_2,\pi _2,\sigma
_2\right) $ be two complex doubles constructed previously such that there
are isomorphisms $f_i:\left( Y,\pi ,\sigma \right) \rightarrow \left(
X_i,\pi _i,\sigma _i\right) $ between the double covers. In both cases its
lifting is $\left( \widehat{X_c},\widehat{\pi _c},\widehat{\sigma _c}\right) 
$ and we denote the projections by $p_i:\widehat{X_c}\rightarrow X_i$, for $%
i=1,2$. With these conditions there exists a unique $\widehat{f}:\left( 
\widehat{Y},\widehat{\pi },\widehat{\sigma }\right) \rightarrow \left( 
\widehat{X_c},\widehat{\pi _c},\widehat{\sigma _c}\right) $ isomorphism
between double covers such that $\widehat{f_1}=\widehat{f_2}=\widehat{f}$.
Let $z\in N\left( X,1\right) $ be such that $p_x^{-1}\left( z\right)
=\left\{ z_1,z_2\right\} $, $\widehat{\pi _c}^{-1}\left( z_i\right) =\left\{
w_{i,1},w_{i,2}\right\} $, $p_1\left( w_{1,1}\right) =p_1\left(
w_{2,2}\right) $, $p_1\left( w_{1,2}\right) =p_1\left( w_{2,1}\right) $, $%
p_2\left( w_{1,1}\right) =p_2\left( w_{2,1}\right) $ and $p_2\left(
w_{1,2}\right) =p_2\left( w_{2,2}\right) $. In this case: 
\[
f_1\circ p_y\circ \widehat{f}^{-1}\left( w_{1,1}\right) =p_1\left(
w_{1,1}\right) =p_1\left( w_{2,2}\right) =f_1\circ p_y\circ \widehat{f}%
^{-1}\left( w_{2,2}\right) , 
\]
\[
f_2\circ p_y\circ \widehat{f}^{-1}\left( w_{1,1}\right) =p_2\left(
w_{1,1}\right) =p_2\left( w_{2,1}\right) =f_2\circ p_y\circ \widehat{f}%
^{-1}\left( w_{2,1}\right) , 
\]
and then $p_y\circ \widehat{f}^{-1}\left( w_{1,1}\right) =p_y\circ \widehat{f%
}^{-1}\left( w_{2,2}\right) =p_y\circ \widehat{f}^{-1}\left( w_{2,1}\right) $%
. But this is impossible because $\#p_y^{-1}\left( x\right) \leq 2$ and $%
\widehat{f}$ is a bijection. Hence we have the unicity of $\left( X_c,\pi
_c,\sigma _c\right) $. The unicity of $f$ is a consequence of the unicity of 
$\widehat{f}$.\TeXButton{End Proof}{\endproof  }

We conclude this section showing diferent situations of complex doubles:

\begin{example}
Let be $\widehat{X}=\widehat{{\Bbb C}^{+}}\times \left\{ 1\right\} \cup 
\widehat{{\Bbb C}^{+}}\times \left\{ 2\right\} $, $X=\widehat{X}/\sim $
where $\sim $ are the identifications: $\left( 0,1\right) \sim \left(
0,2\right) $ and $\left( i,1\right) \sim \left( i,2\right) $, $\widehat{X_c}=%
\widehat{{\Bbb C}}\times \left\{ 1\right\} \cup \widehat{{\Bbb C}}\times
\left\{ 2\right\} $. We define: 
\[
\begin{array}{cccc}
\widehat{\sigma _c}: & \widehat{X_c} & \rightarrow & \widehat{X_c} \\ 
& \left( z,n\right) & \rightarrow & \left( \overline{z},n\right)
\end{array}
,\qquad \qquad 
\begin{array}{cccc}
\widehat{\pi _c}: & \widehat{X_c} & \rightarrow & \widehat{X} \\ 
& \left( z,n\right) & \rightarrow & \left( \phi \left( z\right) ,n\right)
\end{array}
\]
where $\phi \left( z\right) =\func{Re}z+\left| \func{Im}z\right| i$. Let be
now $X_1=\widehat{X_c}/\sim _1$ where $\sim _1$ are the identifications: 
\[
\left( 0,1\right) \sim _1\left( 0,2\right) ,\qquad \left( i,1\right) \sim
_1\left( i,2\right) ,\qquad \left( -i,1\right) \sim _1\left( -i,2\right) , 
\]
and $X_2=\widehat{X_c}/\sim _2$ where $\sim _2$ are the identifications: 
\[
\left( 0,1\right) \sim _2\left( 0,2\right) ,\qquad \left( i,1\right) \sim
_2\left( -i,2\right) ,\qquad \left( -i,1\right) \sim _2\left( i,2\right) . 
\]
Finally we define $\sigma _i$ and $\pi _i$ in the obvious way. Hence $\left(
X_1,\pi _1,\sigma _1\right) $ and $\left( X_2,\pi _2,\sigma _2\right) $ are
the two complex doubles of $X$. Let 
\[
\begin{array}{cccc}
f: & \left( X_1,\sigma _1\right) & \rightarrow & \left( X_2,\sigma _2\right)
\\ 
& \left( z,1\right) & \rightarrow & \left( z,1\right) \\ 
& \left( z,2\right) & \rightarrow & \left( -z,2\right) ,
\end{array}
\]
be an isomorphism between symmetric Riemann surfaces with nodes but it is
not an isomorphism between double covers because $\pi _2\circ f\left(
z,2\right) \neq \pi _1\left( z,2\right) $. In this case the two complex
doubles are isomorphic as symmetric Riemann surfaces with nodes.
\end{example}

\begin{example}
Let $\widehat{X}$ be $\widehat{{\Bbb C}^{+}}$, $X$ be $\widehat{X}/\sim $
where $\sim $ are the identifications: $i\sim 2i$, $\widehat{X_c}$ be $%
\widehat{{\Bbb C}}$ and we define: 
\[
\begin{array}{cccc}
\widehat{\sigma _c}: & \widehat{X_c} & \rightarrow & \widehat{X_c} \\ 
& z & \rightarrow & \overline{z}
\end{array}
,\qquad \qquad 
\begin{array}{cccc}
\widehat{\pi _c}: & \widehat{X_c} & \rightarrow & \widehat{X} \\ 
& z & \rightarrow & \phi \left( z\right) .
\end{array}
\]
Let $X_1$ be $\widehat{X_c}/\sim _1$ where $\sim _1$ are the
identifications: 
\[
i\sim _12i,\qquad -i\sim _1-2i, 
\]
and $X_2$ be $\widehat{X_c}/\sim _2$ where $\sim _2$ are the
identifications: 
\[
i\sim _2-2i,\qquad -i\sim _22i. 
\]
Finally we define $\sigma _i$ and $\pi _i$ in the obvious way. Then $\left(
X_1,\pi _1,\sigma _1\right) $ and $\left( X_2,\pi _2,\sigma _2\right) $ are
the two complex doubles of $X$. It is easy to see that the two complex
doubles are homeomorphic. Since $Fix\sigma _1=Fix\sigma _2={\Bbb R}\cup
\left\{ \infty \right\} $ then $\left( X_1,\sigma _1\right) $ and $\left(
X_2,\sigma _2\right) $ are not homeomorphic as symmetric Riemann surfaces
with nodes because $X_1\backslash Fix\sigma _1$ is not connected and $%
X_2\backslash Fix\sigma _2$ is connected. In this case the two complex
doubles are homeomorphic as Riemann surfaces with nodes but are not
homeomorphic as symmetric Riemann surfaces with nodes.
\end{example}

\begin{example}
Let $\widehat{X}$ be $\widehat{{\Bbb C}}$, $X$ be $\widehat{X}/\sim $ where $%
\sim $ is the identification: $1\sim 2$ and $0\in N\left( X,2\right) $, $%
\widehat{X_c}$ be $\widehat{{\Bbb C}}\times \left\{ 1\right\} \cup \widehat{%
{\Bbb C}}\times \left\{ 2\right\} $ and we define: 
\[
\begin{array}{cccc}
\widehat{\sigma _c}: & \widehat{X_c} & \rightarrow & \widehat{X_c} \\ 
& \left( z,1\right) & \rightarrow & \left( \overline{z},2\right) \\ 
& \left( z,2\right) & \rightarrow & \left( \overline{z},1\right)
\end{array}
,\qquad \qquad 
\begin{array}{cccc}
\widehat{\pi _c}: & \widehat{X_c} & \rightarrow & \widehat{X} \\ 
& \left( z,1\right) & \rightarrow & z \\ 
& \left( z,2\right) & \rightarrow & \overline{z}.
\end{array}
\]
Let $X_1$ be $\widehat{X_c}/\sim _1$ where $\sim _1$ are the
identifications: 
\[
\left( 0,1\right) \sim _1\left( 0,2\right) ,\qquad \left( 1,1\right) \sim
_1\left( 2,1\right) ,\qquad \left( 1,2\right) \sim _1\left( 2,2\right) , 
\]
and $X_2$ be $\widehat{X_c}/\sim _2$ where $\sim _2$ are the
identifications: 
\[
\left( 0,1\right) \sim _2\left( 0,2\right) ,\qquad \left( 1,1\right) \sim
_2\left( 2,2\right) ,\qquad \left( 1,2\right) \sim _2\left( 2,1\right) . 
\]
Finally we define $\sigma _i$ and $\pi _i$ in the obvious way. Then $\left(
X_1,\pi _1,\sigma _1\right) $ and $\left( X_2,\pi _2,\sigma _2\right) $ are
the two complex doubles of $X$ that are not homeomorphic.
\end{example}

The reader can observe that we have not defined when a Klein surface with
nodes is stable. We shall do it now. Let $X$ be a Klein surface with nodes
and $\left( X_1,\pi _1,\sigma _1\right) $, $\left( X_2,\pi _2,\sigma
_2\right) $ be two complex doubles of $X$, then $\chi \left( X_1\right)
=\chi \left( X_2\right) $ and $\#N\left( X_1\right) =\#N\left( X_2\right) $.
Since $X_1$ and $X_2$ are orientable then $g\left( X_1\right) =g\left(
X_2\right) $ and we define the algebraic genus of $X$ by the genus of anyone
of its complex doubles.

Finally, as $\widehat{X_1}=\widehat{X_2}$, then the parts of $X_1$ are
homeomorphic to the parts of $X_2$ and $\left( X_1,\sigma _1\right) $ is a
stable symmetric Riemann surface if and only if $\left( X_2,\sigma _2\right) 
$ is. These results gives the possibility for the following definition:

\begin{definition}
Let $X$ be a Klein surface with nodes. We say that $X$ is a stable Klein
surface if anyone of its complex doubles is a stable symmetric Riemann
surface.
\end{definition}


\begin{thebibliography}{A-G}
\bibitem[A-G]{A-G}  N.L. Alling, N. Greenleaf, {\it Foundations of the
Theory of Klein Surfaces}. Lecture Notes in Math., 219, Springer, 1971.

\bibitem[D-M]{D-M}  P. Deligne, D. Mumford, {\it The irreducibility of the
space of curves of given genus}. Inst. Hautes \'{E}tudes Sci. Publ. Math., 
{\bf 36}, (1969), 75-109.

\bibitem[G]{G}  I. Garijo, {\it Superficies de Riemann y Klein con nodos}.
Ph. D. Thesis, U.N.E.D., 2000.

\bibitem[K]{K}  M. Kontsevich, {\it Intersection theory on the moduli space
of curves and the matrix Airy function}. Communications in Mathematical
Physics, {\bf 147}, (1992), 1-23.

\bibitem[L-Z]{L-Z}  S. K. Lando and A. K. Zvonkin, {\it Graphs on surfaces
and their applications}. Encyclopedia of Mathematical Sciences, vol. 141,
Springer, Berlin 2004.

\bibitem[Se]{Se}  M. Sepp\"{a}la, {\it Moduli space of stable real algebraic
curves}. Ann. Scient. \'{E}c. Norm. Sup., 4$^e$ s\'{e}rie, {\bf 24}, (1991),
519-544.

\bibitem[Si]{Si}  R. Silhol, {\it Compactification of moduli spaces in real
algebraic geometry}. Invent. Math., {\bf 107}, (1992), 151-202.

\bibitem[Z]{Z}  D. Zvonkine, {\it Strebel differentials on stable curves and
Kontsevich's proof of Witten's conjeture}, Preprint, arXiv:\ math.AG/0209071
v2 7 Jan 2004.
\end{thebibliography}
\end{document}